\documentclass[12pt,a4paper,draft]{amsart}
\usepackage{epsf,amsfonts,amscd,latexsym}
\usepackage[T2A]{fontenc}
\usepackage[english]{babel}
\usepackage[matrix,arrow]{xy}

\title[Birational rigidity is not an open property]
{Birational rigidity is not an open property}
\thanks{This work was supported by the grants RFBR 05-01-00353a, NSh-9969.2006.1, MD-4261.2006.1, CRDF RUM1-2692-MO-05}
\author{Ivan Cheltsov}
\email{chelstov@yahoo.com}
\author{Mikhail Grinenko}
\email{grin@mi.ras.ru}
\date{}

\newtheorem{theorem}{\sc Theorem}[section]
\newtheorem{proposition}[theorem]{\sc Proposition}
\newtheorem{lemma}[theorem]{\sc Lemma}
\newtheorem{corollary}[theorem]{\sc Corollary}

\newtheorem{definition}[theorem]{\sc Definition}

\makeatletter

\@addtoreset{equation}{section}
\newcommand{\l@abcd}[2]{\hbox to\textwidth{#1\dotfill #2}}

\makeatother

\newcommand*{\mybegintheorem}[1]{\begin{trivlist}\it%
  \item[\hspace{\labelsep}{\bf #1}]}
  \newcommand*{\myendtheorem}{\end{trivlist}}
  \newenvironment*{theorem*}{\mybegintheorem{Theorem.}}{\myendtheorem}
  \newenvironment*{proposition*}{\mybegintheorem{Proposition.}}{\myendtheorem}
  \newenvironment*{corollary*}{\mybegintheorem{Corollary.}}{\myendtheorem}
  \newenvironment*{definition*}{\mybegintheorem{Definition.}}{\myendtheorem}

\theoremstyle{remark}
\newtheorem{remark}[theorem]{\sc Remark}

\renewcommand{\phi}{\varphi}
\renewcommand{\epsilon}{\varepsilon}

\newcommand{\lra}{\longrightarrow}

\newcommand{\PQ}{{{\mathbb P}^4}}
\newcommand{\PT}{{{\mathbb P}^3}}
\newcommand{\PTw}{{{\mathbb P}^2}}
\newcommand{\POn}{{{\mathbb P}^1}}
\newcommand{\ZA}{{\mathbb Z}}
\newcommand{\QA}{{\mathbb Q}}

\newcommand{\CA}{{\mathbb C}}
\newcommand{\PA}{{\mathbb P}}
\newcommand{\FA}{{\mathbb F}}
\newcommand{\mc}{\mathcal}
\newcommand{\mf}{\mathfrak}

\newcommand{\Supp}{\mathop{\rm Supp}\nolimits}
\newcommand{\mult}{\mathop{\rm mult}\nolimits}

\newcommand{\rk}{\mathop{\rm rk}\nolimits}
\newcommand{\Pic}{\mathop{\rm Pic}\nolimits}
\newcommand{\Cl}{\mathop{\rm Cl}\nolimits}
\newcommand{\Bir}{\mathop{\rm Bir}\nolimits}

\newcommand{\eps}{\varepsilon}

\begin{document}
\begin{abstract}
We construct an example of the birationally rigid complete intersection of a quadric and a cubic in $\PA^5$ with an ordinary double point, which under a small deformation gives a non-rigid Fano variety. Thus we show that birational rigidity is not open in moduli.
\end{abstract}

\maketitle

\section{Introduction.}
\label{pre_sec}

In this article we discuss the question, which is closely related to the nature of the birational rigidity notion: whether birational rigidity is open in moduli, or not?

We consider all varieties to be defined over an algebraically closed field of characteristic 0 (e.g., over $\CA$). We recall that a triple $\mu:V\to S$ is said to be a {\it Mori fibration} (also Mori fiber space) if $V$ is a projective $\QA$-factorial terminal threefold, $S$ is a projective normal variety with $\dim S<\dim X$, and $\mu$ is an extremal contraction of fibering type, i.e., the relative Picard number $\rho(V/S)=\rk\Pic(V)-\rk\Pic(S)$ is equal to 1 and $(-K_V)$ is $\mu$-ample. Abusing the notation, we will denote Mori fibrations also $V\to S$, $V/S$, or simply $V$ when the corresponding contractions or bases are clear. We have the following possibilities for Mori fibrations: Fano varieties ($\dim S=0$), del Pezzo fibrations ($\dim S=1$), and conic bundles ($\dim S=2$).

We say that a birational map $\chi:X\dasharrow X'$ of two Mori fibrations $X\to S$ and $X'\to S$ is {\it square} if it fits into a commutative diagram
$$
\begin{array}{ccc}
X & \stackrel{\chi}{\dasharrow} & X' \\
\downarrow && \downarrow \\
S & \stackrel{\psi}{\lra} & S'
\end{array}
$$
where $\psi$ is a birational map, and moreover, $\chi$ induces an isomorphism $X_{\eta}\stackrel{\simeq}{\lra}X'_{\eta}$ of the fibers over the generic point $\eta$ (or, which is the same, $\chi$ induces isomorphisms of general fibers).

\begin{definition}
\label{rig_def}
A Mori fibration $X/S$ is said to be birationally rigid, if for any birational map $\phi:X\dasharrow X'$ to another Mori fibration $X'/S'$ there is a birational map $\mu\in\Bir(X)$ such that the composition $\phi\circ\mu:X\dasharrow X'$ is square, and $X/S$ is birationally superrigid, if $\mu$ can be chosen biregular (or simply an identity map).
\end{definition}

Consider a Fano variety $X$ and a birational map $\phi:X\dasharrow X'$ to a Mori fibration $X'/S'$. If $X$ is (birationally) rigid, the definition says that $X'$ is isomorphic to $X$ and $\phi$ can be viewed as a birational automorphism. If $X$ is superrigid, then $\phi$ is an isomorphism itself.

One of the most common conjecture about birationally rigid varieties was that small deformations keep the rigidity (e.g., \cite{Co}, conjecture 1.4): given any scheme $T$, and a flat family of Mori fibrations ${\mc X}\to{\mc S}$ parameterized by $T$, the set of all $t\in T$ such that the corresponding fiber ${\mc X}_t\to{\mc S}_t$ is birationally rigid, is open in $T$ (possibly empty). In other words, birational rigidity is open in moduli. Up to this moment, all known examples of rigid varieties satisfy this conjecture.

Nevertheless, in this paper we show that the conjecture falls. Our counter-examples are based on degenerations of the complete intersections of quadrics and cubics in $\PA^5$. It is known (\cite{Isk}, \cite{IP}) that a general non-singular Fano variety of this kind is rigid. We construct a special class of the singular complete intersections that are generally non-rigid, but it contains a sub-class of rigid varieties. Our consideration is based on ideas of V.A.Iskovskikh and A.V.Pukhlikov's works
, and we put only those places of the proof that are different from their explanations. At the last section we give counter-examples in the class of del Pezzo fibrations.

The essential part of the work was done during the stay of one of the author at the Max-Planck-Institute f\"ur Mathematik in August-September 2006, and we would like to express our gratitude to the Institute for hospitality and excellent condition for work. The authors also thank K.Shramov for useful conversations.

\section{Main result.}
\label{main_sec}

We construct the essential family of the singular complete intersections of quadrics and cubics as follows. Let $X\subset\PQ$ be a quartic given by the equation
\begin{equation}
\label{X_eq}
ht-q_1q_2=0,
\end{equation}
where $h$, $t$, $q_1$, and $q_2$ are homogeneous polynomials of degree 1, 3, 2, and 2 respectively, in the coordinates $[y_0:y_1:\ldots:y_4]$ in $\PQ$. We assume that $X$ has the only 12 different ordinary double points given by the equalities
$$
h=t=q_1=q_2=0.
$$
It is easy to see that the hyperplane $\{h=0\}$ cuts off two quadratic surfaces $S_1=\{h=q_1=0\}$ and $S_2=\{h=q_2=0\}$. We have
$$
\begin{array}{lll}
\Pic(V)&=&\ZA[-K_X], \\
\Cl(V)&=&\ZA[-K_X]\oplus\ZA[S_1],
\end{array}
$$
thus $X$ is not $\QA$-factorial, so it does not belong to a category of Mori fibrations. We can get two Mori fibrations that are birational to $X$, using two "unprojections" as follows. Let $[y_0:y_1:\ldots:y_5]$ be the coordinates in $\PA^5$, and consider varieties $V_1$ and $V_2$ defined by the equations
\begin{equation}
\label{V1V2_eq}
\begin{array}{cc}
V_1= \left\{
\begin{array}{ccl}
y_5h & = & q_1 \\
y_5q_2 & = & t
\end{array}
\right. , & V_2= \left\{
\begin{array}{ccl}
y_5h & = & q_2 \\
y_5q_1 & = & t
\end{array}
\right.
\end{array}
\end{equation}
Note that the polynomials $h$, $t$, $q_1$, and $q_2$ do not depend on $y_5$, and we get the equation (\ref{X_eq}) by excluding the variable $y_5$ from the equations for $V_1$ or $V_2$. In other words, we can get $X$ from $V_1$ or $V_2$ by the projection $\PA^5\dasharrow\PA^4$ from the point $\xi=(0:\ldots:0:1))$.

It is easy to check that under a general choice of the corresponding polynomials the varieties $V_1$ and $V_2$ are Fano varieties with a unique ordinary double point $\xi$. The birational maps $\phi_1:V_1\dasharrow X$, $\phi_2:V_2\dasharrow X$, and $\psi:V_1\dasharrow V_2$ can be described as follows. Let $\alpha_1:\tilde V_1\to X$ be the blow-up of the quadric $S_2$. Thus $\alpha_1$ is the small resolution of singularities of $X$ with the exceptional lines $\tilde l_1,\ldots,\tilde l_{12}$. Denote $\tilde S_1$ and $\tilde S_2$ the strict transforms of the corresponding quadratic surfaces. The variety $\tilde V_1$ is non-singular, and the quadric $\tilde S_2$ has the normal sheaf ${\mc O}(-1)$. So there exists the divisorial contraction $\beta_1:\tilde V_1\to V_1$ of the exceptional divisor $\tilde S_2$. Now we set $\phi_1=\alpha_1\circ\beta_1^{-1}$. The birational map $\phi_2$ is constructed in the same way as $\phi_1$, but first we blow up the quadric $S_2$. Finally, there 
 exists a flop $\gamma:\tilde V_1\dasharrow\tilde V_2$ centered at all 12 lines $\tilde l_1,\ldots,\tilde l_{12}$ simultaneously, and we have $\psi=\beta_2\circ\gamma\circ\beta_1^{-1}: V_1\dasharrow V_2$.

Denote ${\mc F}$ a family of all the complete intersections of quadrics and cubics in $\PA^5$ that are constructed like $V_1$ and $V_2$. So we see that each variety $U_1\in{\mc F}$ has its birational "counterpart" $U_2$ like $V_1$ and $V_2$, and $U_1$ is not isomorphic to $U_2$ in general. Indeed, any isomorphism of $U_1$ and $U_2$ is induced by an automorphism of $\PA^5$ that keep fixed the singular point $(0:\ldots:0:1)$. Thus the projection from this point induces also an involution of $\PQ$. The equation of the corresponding quartic is invariant with respect to the involution, which change the places of $q_1$ and $q_2$ in (\ref{X_eq}). It is impossible in general case (e.g., we can choose an equation for $h$ which is not invariant with respect to any such an involution).

Now we construct the subfamily ${\mc F}_r\subset{\mc F}$. Consider a reflection $\iota:\PQ\to\PQ$ given by the following action:
\begin{equation}
\label{Xinv_eq}
\left\{
\begin{array}{ll}
y_0\to -y_0 & \\
y_i\to y_i, & i=1\ldots 4
\end{array}
\right.
\end{equation}
We continue the action of $\iota$ assuming $y_5\to y_5$. Let $h$ and $t$ be homogeneous polynomials of degree 1 and 3 in the variables $y_1,\ldots,y_4$, and $q$ be a homogeneous polynomial of degree 2 in the variables $y_0,\ldots,y_4$. Thus $h$ and $t$ are invariant with respect to the reflection $\iota$. Consider a quadric $X\subset\PQ$ given by the equation
$$
ht-q\iota^*(q)=0.
$$
We can choose the birational modifications of $X$ as before, and get Fano varieties
\begin{equation}
\label{V1V2inv_eq}
\begin{array}{cc}
V_1= \left\{
\begin{array}{ccl}
y_5h & = & q \\
y_5\iota^*(q) & = & t
\end{array}
\right. , & V_2= \left\{
\begin{array}{ccl}
y_5h & = & \iota^*(q) \\
y_5q & = & t
\end{array}
\right.
\end{array}
\end{equation}
We see immediately that $\iota(V_1)=V_2$, i.e., $V_1$ and $V_2$ are isomorphic, and the birational map $\psi:V_1\dasharrow V_2$ is actually a birational automorphism.

Let us fix the reflection $\iota$, and suppose ${\mc F}_r$ consisting of varieties of the kind (\ref{V1V2inv_eq}). Obviously, ${\mc F}_r\subset{\mc F}$, and we can make an important observation that varieties from ${\mc F}$ can be viewed as small deformations of varieties from ${\mc F}_r$.

Now we describe some conditions of generality for varieties of these two families. Let $Q$ and $T$ be a quadric and a cubic in $\PA^5$, and $V$ their complete intersection with a unique ordinary double point $\xi$. We assume that $V$ satisfies the following additional conditions of generality:
\begin{itemize}
\item[${\mc G}_1$]: the quadric $Q$ is non-singular;%
\item[${\mc G}_2$]: if $l\subset V$ is a line and $P\in\PA^5$ is a plane that contains $l$, then the scheme-theoretical intersection $V\cap P$ is reduced along $l$ (in case $\xi\not\in l$ it means that the normal sheaf ${\mc N}_{l|V}\simeq{\mc O}\oplus{\mc O}(-1)$, see \cite{IP}, chapter 3, proposition 1.1);%
\item[${\mc G}_3$]: for any plane $P\subset\PA^5$, the intersection $V\cap P$ is not three lines with a common point, and if $\xi\in P$, it does not consist of any three lines;%
\item[${\mc G}_4$]: given $l\subset L\subset\PA^5$, where $l\subset V$ is a line through the point $\xi$ and $L$ is a three-dimensional subspace such that $Q|_L$ consists of two planes (with the common line $l$), then $L$ is not the tangent space to $V$ for any point in $l\setminus\{\xi\}$;%
\item[${\mc G}_5$]: let $l\subset V$ be a line, $\xi\in l$, then for each point $B\in l\setminus\{\xi\}$ there are not more than 3 lines passing through $B$ (including $l$ itself).%
\end{itemize}

\noindent Our main result is the following:
\begin{theorem}
\label{main_th}
(i) A general variety from ${\mc F}_r$ or ${\mc F}$ satisfies the conditions ${\mc G}_1,\ldots,{\mc G}_5$.

\noindent(ii) Let $V_1$ and $V_2$ be general varieties from ${\mc F}$ or ${\mc F}_r$ given by the equations (\ref{V1V2_eq}). Then $V_1$ and $V_2$ are unique Mori fibrations in their class of birational equivalence, and they are isomorphic if they belong to the subfamily ${\mc F}_r$.
\end{theorem}
\noindent From this theorem we can deduce immediately:
\begin{corollary}
\label{main_cor}
Birational rigidity is not open in moduli.
\end{corollary}
Indeed, by the part (ii) of theorem \ref{main_th}, general varieties from the family ${\mc F}_r$ are birationally rigid, but general varieties from ${\mc F}$ are non-rigid. So we conclude that birational rigidity is not kept by small deformations.

Thus, it is proved that a general $V\in{\mc F}$ has exactly 2 different models of Mori fibration, both of them are Fano varieties. The first example of this kind was constructed by A.Corti and M.Mella in \cite{CoMe}. From the viewpoint of their work, general varieties from ${\mc F}$ have the pliability to be equal to 2, and equal to 1 for varieties from ${\mc F}_r$.

The further exposition is organized as follows. We prove the part (i) of theorem \ref{main_th} in section \ref{gencond_sec}, and the part (ii) in sections \ref{excl1_sec}, \ref{excl2_sec}, and \ref{fin_sec}. We conclude with some further examples in section \ref{rem_sec}.

\section{Generality conditions.}
\label{gencond_sec}

In this section we prove the part (i) of theorem (\ref{main_th}) in lemmas \ref{g1cond_lem} and \ref{g2-g5cond_lem}. Is is clear that it is enough to prove the generality conditions ${\mc G}_1,\ldots,{\mc G}_5$ only for varieties from ${\mc F}_r$, the result for ${\mc F}$ will follow automatically.

Let $[y_0:y_1:\ldots:y_5]$ be the coordinates in $\PA^5$. We assume the reflection $\iota$ to be defined by
$$
\left\{
\begin{array}{ll}
y_0\to -y_0 & \\
y_i\to y_i, & i=1,\ldots,5.
\end{array}
\right.
$$
The set of fixed points of $\iota$ are the plane $F=\{y_0=0\}$ and the point $\zeta=(1:0:\ldots:0)$. We choose the polynomials $h, t, q$ as follows:
$$
\begin{array}{rcl}
h(y_*)&=&\sum_{i=1}^4\alpha_iy_i, \\
q(y_*)&=&Ay_0^2+y_0\sum_{i=1}^4\beta_iy_i+\sum_{1\le i\le j\le 4}\gamma_{ij}y_iy_j, \\
t(y_*)&=&\sum_{1\le i\le j\le k\le 4}\delta_{ijk}y_iy_jy_k+
y_0^2\sum_{i=1}^4\eps_iy_i.
\end{array}
$$
Assuming a variety $V\in{\mc F}_r$ to be the complete intersection of a quadric $Q$ and a cubic $T$, we have the following equations for these varieties:
\begin{equation}
\label{yV_eq}
\left\{
\begin{array}{l}
y_5\sum\alpha_iy_i=Ay_0^2+y_0\sum\beta_iy_i+\sum\gamma_{ij}y_iy_j; \\
y_5\left(Ay_0^2-y_0\sum\beta_iy_i+\sum\gamma_{ij}y_iy_j\right)=
\sum\delta_{ijk}y_iy_jy_k+y_0^2\sum\eps_iy_i;
\end{array}
\right.
\end{equation}
where the coefficients $A$, $\alpha_i$, $\beta_i$, $\gamma_{ij}$, $\delta_{ijk}$, and $\eps_i$ are the parameters.

Notice that $V$ has a singular point at $\xi=(0:\ldots:0:1)$. We will always assume that the coefficient $A$ is not 0, thus $V$ does not pass through the point $\zeta$ (the isolated fixed point of the reflection $\iota$).

\begin{lemma}
\label{cond1_lem}
Let $l\subset V$ be a line that does not pass trough the singular point $\xi\in V$. Then $l$ does not intersect the line $L\subset\PA^5$ through the points $\zeta$ and $\xi$.
\end{lemma}
\noindent{\bf Proof.} Assume the converse, i.e., $l\cap L\ne\emptyset$. Then there exists a plane $P$ such that $l,L\subset P$. We may suppose that $l$ intersects the hyperplane $F$ at the point $(0:1:0:0:0:0)$, and since $L=\{y_1=y_2=y_3=y_4=0\}$, the plane $P$ is defined by $\{y_2=y_3=y_4=0\}$. Notice that $A\ne 0$ by the assumption, so $P\not\subset Q$.

The restriction $Q|_P$ gives us a conic
$$
\alpha_1y_1y_5=Ay_0^2+\beta_1y_0y_1+\gamma_{11}y_1^2,
$$
as it follows from the equation (\ref{yV_eq}). On the other hand, $l\subset Q|_P$, so the equation before have to be a production of two linear forms. From this we deduce that $\alpha_1=0$, thus $\xi\in l$, which contradicts to the conditions. Lemma \ref{cond1_lem} is proved.

\begin{lemma}
\label{cond2_lem}
Let $V$ be a general variety in ${\mc F}_r$, and $l\subset V$ a line that does not pass through the singular point $\xi\in V$. Then $l\not\subset F$.
\end{lemma}
\noindent{\bf Proof.} We will argue by counting of dimensions. Since $\xi\not\in l$, we can take the projection $\PA^5\dasharrow \PQ$ and prove the lemma the the quartic $X$, which is the image of $V$.

Suppose $l'\subset F'$, where the line $l'$ and the hyperplane $F'=\{y_0=0\}$  are the images of $l$ and $F$ respectively. Then the restriction $X'=X|_{F'}$ has the equation
$$
\left(\sum\alpha_iy_i\right)\left(\sum\delta_{ijk}y_iy_jy_k\right)+
\left(\sum\gamma_{ij}y_iy_j\right)^2=0.
$$
We may suppose that $l'$ is defined by the equations $y_3=y_4=0$ in $F'$, so from $l'\subset X$ we have the following conditions:
$$
\begin{array}{l}
\alpha_1\delta_{111}+\gamma_{11}^2=0, \\
\alpha_1\delta_{112}+\alpha_2\delta_{111}+2\gamma_{11}\gamma_{12}=0, \\
\alpha_1\delta_{122}+\alpha_2\delta_{112}+2\gamma_{11}\gamma_{22}+\gamma_{12}^2=0, \\
\alpha_1\delta_{222}+\alpha_2\delta_{122}+2\gamma_{12}\gamma_{22}=0, \\
\alpha_2\delta_{222}+\gamma_{22}^2=0.
\end{array}
$$
It is not very difficult to check that we they gives 5 independent conditions. Let ${\mc X}'$ be the set of all two-dimensional quartics that are defined by equations of the same form as $X'$. Consider the production $G(2,4)\times{\mc X}'$ with the corresponding projections $p$ and $q$ onto $G(2,4)$ and ${\mc X}'$ respectively, and a subvariety $I=\{(l',X'): l'\subset X'\}\subset G(2,4)\times{\mc X}'$. Now we have
$$
\dim I=\dim{\mc X'}-5+\dim G(2,5)=\dim{\mc X}'-1 < \dim{\mc X}'.
$$
Thus a general variety $X$, and hence a general $V\in{\mc F}_r$, has no lines with the indicated conditions. Lemma \ref{cond2_lem} is proved.

\begin{lemma}
\label{g1cond_lem}
A general variety $V\in{\mc F}_r$ satisfies condition ${\mc G}_1$.
\end{lemma}
\noindent{\bf Proof.} The lemma is obvious: the quadric $Q$ is non-singular a for general choice of the parameters in the equation (\ref{yV_eq}).

In what follows, instead of fixing the reflection $\iota$, we will fix the hyperplane $F$, the point $\xi$, which will be the singular point for our varieties, and the line $l\not\subset F$ (see lemma \ref{cond2_lem}), and then will move in $\PA^5$ the fixed point $\zeta$ of the reflection. Notice also that any linear transformation of the coordinates $y_1,\ldots,y_4$ keeps the look of the equations (\ref{yV_eq}).

Thus we fix a new coordinates $[x_0:\ldots:x_5]$ in $\PA^5$ and assume $F=\{x_0=0\}$, $l=\{x_2=\ldots=x_5=0\}$, $\xi=(0:\ldots:0:1)$, $\zeta=(1:-a_1:-a_2:-a_3:-a_4:-a_5)$. The relation between the two system of coordinates is given by
\begin{equation}
\label{xycoord_eq}
\left\{
\begin{array}{l}
y_0=x_0;  \\
y_i=x_i+a_ix_0, \quad i=1,\ldots,4; \\
y_5=x_5+bx_1+(a_5+ba_1)x_0;
\end{array}
\right.
\end{equation}
where $b$ and $a_i$ are some numbers. The reflection $\iota$ acts now as follows:
$$
\begin{array}{l}
x_0\to -x_0; \\
x_i\to x_i+2a_ix_0, \quad i=1,\ldots,5.
\end{array}
$$

\begin{lemma}
\label{g2-g5cond_lem}
A general variety $V\in{\mc F}_r$ satisfies conditions ${\mc G}_2$, ${\mc G}_3$, ${\mc G}_4$, and ${\mc G}_5$.
\end{lemma}

\noindent{\bf Proof.} In all cases, we argue by counting the dimensions. First we consider condition ${\mc G}_2$, the case $\xi\not\in l$. Let $l\subset V$ be a line, and $L\subset\PA^5$ be a line that passes trough the points $\xi$ and $\zeta$. By lemmas \ref{cond1_lem} and \ref{cond2_lem}, for any $l$ we may assume $l\not\subset F$ and $l\cap L=\emptyset$. Consider an open subset of a variety of $(1,2)$--flags $T=\{(l,P): l\subset P, l\not\subset F, l\cap L=\emptyset\}$, and a closed subset $S=\{(l,P): P\subset\langle l, L\rangle\}\subset T$, where $\langle l, L\rangle$ denote a unique 3-dimensional linear subspace that contains both the lines. It easy to compute that $\dim T=11$ and $\dim S=6$.

Consider a flag $(l,P)\in T\setminus S$. Then we choose the coordinates $[x_0:\ldots:x_5]$ as before this lemma, and we can assume $l=\{x_2=x_3=x_4=x_5=0\}$ and $P=\{x_3=x_4=x_5=0\}$. It is important to observe that the numbers $a_3$ and $a_4$ can not vanish simultaneously because $P\not\subset\langle l,L\rangle$. Suppose $l\subset V$ and $V\cap P$ is not reduced along $l$. Then we substitute the coordinates $[y_*]$ for $[x_*]$ in the equation (\ref{yV_eq}) using the system (\ref{xycoord_eq}), and look at the restrictions $Q|_P$ and $T|_P$. By the assumption, these restrictions are not reduced along $l$. It is not very difficult to see that this gives 12 independent linear conditions for the coefficients in (\ref{yV_eq}). Consider a subvariety $I\subset {T\setminus S}\times{\mc F}_r$ consisting of all pairs $\left((l,P),V\right)$ such that $l\subset V$ and $V\cap P$ is not reduced along $l$. Then $\dim I=\dim T\setminus S+\dim{\mc F}_r-12=\dim{\mc F}_r-1$, thus $I$ can not co
 ver ${\mc F}_r$ under the projection ${T\setminus S}\times{\mc F}_r\to{\mc F}_r$. We argue by the same way for the situation $(l,P)\in S$, with a few modifications. Exactly, here $a_3=a_4=0$, but $a_2\ne 0$ since $l\cap L=\emptyset$. Depending on the mutual location of $\xi$ and $P$, we obtain 9 or 10 linear conditions for the coefficients in (\ref{yV_eq}). It remains to take into account that $\dim S=6$. So finally we see that the condition holds for a general $V$ and lines $l$ such that $\xi\not\in l$.

The case $\xi\in l\subset V$ for the condition ${\mc G}_2$ and the conditions ${\mc G}_3$, ${\mc G}_4$, and ${\mc G}_5$ can be proved in the same way by counting the dimensions. Lemma \ref{g2-g5cond_lem} is proved.

\section{Maximal singularities centered at smooth points.}
\label{excl1_sec}

Consider a general variety $V\in{\mc F}$. We assume that $V$ satisfies the conditions ${\mc G}_1$--${\mc G}_5$. In this section we prove the following result:

\begin{proposition}
\label{point_prop}
Let ${\mc D}\subset|n(-K_V)|$ be a linear system on $V$ without fixed components. Suppose ${\mc D}$ has no maximal singularities centered at curves. Then ${\mc D}$ has no maximal singularities centered at points in $V\setminus\{\xi\}$.
\end{proposition}

We use the method of maximal singularities (\cite{Pukh1}). Suppose that ${\mc D}$ has a maximal singularity centered at the point $B_0\in V\setminus\{\xi\}$. It means that there exists a discrete valuation ${\mf v}$ centered at $B_0$ such that the N\"other-Fano inequality for ${\mc D}$ with respect to ${\mf v}$ holds: ${\mf v}(\mc D)>n\delta_{\mf v}$. Here $\delta_{\mf v}$ means the canonical multiplicity with respect to ${\mf v}$. Then we know (\cite{MKaw1}) that there exists a discrete valuation ${\mf v}_{div}$ that is centered at $B_0$ and can be realized by a weighted blow-up with the weights $(1,L,N)$, where $1\le L<N$ or $N=L=1$. In its turn, the weighted blow-up can be realized as a chain of usual blow-ups
$$
V_N\stackrel{\phi_N}{\lra}V_{N-1}\stackrel{\phi_{N-1}}{\lra}\ldots
\stackrel{\phi_2}{\lra}V_1\stackrel{\phi_1}{\lra}V_0=V
$$
with centers $B_{i-1}\subset V_{i-1}$ and exceptional divisors $E_i\subset V_i$, where:
\begin{itemize}
\item $B_0,\ldots,B_{L-1}$ are points and $B_L,\ldots,B_{N-1}$ are curves;
\item $B_L$ is a line in $E_L\cong\PTw$, $B_i$ for $i>L$ is a section of the corresponding linear surface $E_i\cong\FA_{i+1-L}$ that does not intersect the minimal section;
\item in all cases $B_i\cap E_{i-1}^i=\emptyset$ (upper indices denote the strict transform of a curve or a divisor on the corresponding floor of the chain of blow-ups).
\end{itemize}
Denote $\nu_i=\mult_{B_{i-1}}{\mc D}^{i-1}$, $i>0$. Then the N\"other-Fano inequality looks as follows:
\begin{equation}
\label{NF_ineq}
\nu_1+\ldots+\nu_N>n(L+N).
\end{equation}
We consider the following possible cases:
\begin{itemize}
\item $N>L>1$ -- the general infinitely near case;
\item $N>L=1$ -- the special infinitely near case (we blow up the point $B_0$ and then the curves $B_1,\ldots,B_{N-1}$);
\item $N=L=1$ -- the "point" case (${\mf v}_{div}$ is realized by a single blow-up of the point $B_0$, and then $\nu_1>2n$).
\end{itemize}

Suppose that there is no lines on $V$ that pass through the point $B_0$ or all such lines does not contain the singular point $\xi\in V$. Then we can use the argumentation of the paper \cite{IP}. So in what follows, we assume that $l$ is one of the 12 lines on $V$ that pass through the singular point $\xi$, and $B_0\in l$.

\medskip

\noindent{\bf The general infinitely near case.} Let $D_1$ and $D_2$ be general elements of the linear system ${\mc D}$. Then we can put
$$
D_1\circ D_2=\alpha l+C,
$$
where $C$ is the residual curve, $l\not\subset\Supp C$, and general $H\in|-K_V|$ we have
$$
\deg\left(\alpha l+C\right)=\alpha+\deg C=(\alpha l+C)\circ H=6n^2.
$$
We introduce the number $k$ by
$$
k=\max\left\{i\le L: B_{i-1}\in l^{i-1}\right\}.
$$
Clearly, $k>0$, and either $k=L$ or $B_k\not\in l^k$. Then, we have the so-called quadratic inequality for the cycle $D_1\circ D_2$ (see \cite{Pukh1}):
$$
\sum_{i=1}^L\mult_{B_{i-1}}\left(D_1\circ D_2\right)^{i-1}\ge
\sum_{i=1}^N\nu_i^2>n^2\frac{(N+L)^2}N,
$$
or, denoting $m_i=\mult_{B_{i-1}}C^{i-1}$,
$$
k\alpha+\sum_{i=1}^L m_i >n^2\frac{(N+L)^2}N.
$$
Assume $k=1$, i.e., $B_1\not\in l^1$. Then the linear system $|-K_V-B_0-B_1|$ (those elements from $|-K_V|$ that pass through $B_0$ and the infinitely near point $B_1$) has no basic curves on $V$. For a general $H\in|-K_V-B_0-B_1|$ we have
$$
m_1+m_2\le C\circ H=6n^2-\alpha,
$$
and from the quadratic inequality we get a contradiction:
$$
3n^2L\ge \alpha+(3n^2-\frac12\alpha)L\ge\alpha+\sum_{i=1}^L m_i >
n^2\frac{(N+L)^2}N>4n^2L.
$$
So we assume $k\ge 2$. Using a general element $H$ of the linear system $|-K_V-B_0-\ldots-B_{k-1}|$, we have
$$
m_1+\ldots+m_k\le C\circ H=6n^2-\alpha,
$$
and from the quadratic inequality we get
$$
k\alpha+6n^2-\alpha+\frac{6n^2-\alpha}k(L-k)>n^2\frac{(N+L)^2}N,
$$
or an even more rough estimation,
\begin{equation}
\label{kvad1_ineq}
(k-1)\alpha+6n^2\frac{L}k>n^2\frac{(N+L)^2}N.
\end{equation}

Now one needs to get an upper estimation of $\alpha$. Denote $\mu=\mult_{l}{\mc D}$ and $\nu_0=\mult_{\xi}{\mc D}$. Clearly, $\nu_0\ge\frac12\mu$. Consider the birational morphisms
$$
\tilde V\stackrel{\tilde\phi}{\lra} V'\stackrel{\phi'}{\lra}V_{k-1},
$$
where $\phi'$ is the blow-up of the singular point $\xi$ with the exceptional divisor $E'\cong\FA_0$, and $\tilde\phi$ the blow-up of the strict transform of the line $l$. Denote $\tilde E$ the exceptional divisor of $\tilde\phi$, and $\tilde{\mc H}$ the strict transform of the linear system $|H-B_0-\ldots-B_k-l|$. It is easy to see that $\tilde{\mc H}$ has no basic curves, and moreover, this linear system is ample on $\tilde E\cong\FA_1$ (indeed, $\tilde{\mc H}|_{\tilde E}\subset|s+2f|$, where $s$ and $f$ are the minimal section and a fiber of $\tilde E$ respectively). We observe two important things. First,
\begin{equation}
\label{DE_eq}
\tilde{\mc D}|_{\tilde E}\subset|\mu s+(n+(k+1)-\nu_1-\ldots-\nu_k-\nu_0)f|,
\end{equation}
and denoting $\theta=\frac1k(\nu_1+\ldots+\nu_k)$, we get
\begin{equation}
\label{mu1_ineq}
(k+\frac12)\mu+n\ge k\theta.
\end{equation}
Second, we have
$$
\tilde D_1\circ\tilde D_2\circ\tilde H=
6n^2-k\theta^2-2\mu(n-k\theta)-\nu_0^2-(\mu-\nu_0)^2-(k+1)\mu^2,
$$
hence
$$
\alpha\le\mu_2+\tilde D_1\circ\tilde D_2\circ\tilde H<
6n^2-k\theta^2-2\mu(n-k\theta)-k\mu^2.
$$
This quadratic gets its maximal value at
$$
\mu=\frac{k\theta-n}k,
$$
and we have
\begin{equation}
\label{alpha1_ineq}
\alpha<6n^2-k\theta^2+\frac{(k\theta-n)^2}k.
\end{equation}
Notice that if $\theta\le\frac54n$, from the N\"other-Fano inequality (\ref{NF_ineq}) we obtain $N>4L$, and we get a contradiction with (\ref{kvad1_ineq}). So we assume $\theta>\frac54n$, and taking into account that $\alpha\le 6n^2$, from (\ref{alpha1_ineq}) we have
\begin{equation}
\label{alpha2_ineq}
\alpha<\frac{23}6n^2.
\end{equation}
On the other hand, from (\ref{mu1_ineq}) we obtain
$$
\theta\le\frac{(k+\frac12)\mu+n}{k}\le\frac32n,
$$
and then the N\"other-Fano inequality (\ref{NF_ineq}) yields
$$
N>2L.
$$
This estimation and the inequality (\ref{kvad1_ineq}) give
\begin{equation}
\label{kvad2_ineq}
L\left(\frac{k-1}{L}\alpha+\frac{6n^2}{k}\right)>\frac92n^2L.
\end{equation}
Suppose $L\ge k+1$. Then from (\ref{alpha2_ineq}) and (\ref{kvad2_ineq}) we get $k<2$, a contradiction.

So we assume $L=k\ge 2$. Let us note that there are two lines on the exceptional divisor $E'\cong\FA_0$ that pass through the point $l'\cap E'$, and for at least one of them, say $p$, its strict transform $\tilde p$ on $\tilde V$ does not intersect the minimal section of the divisor $\tilde E$. Denoting by $\eps$ the multiplicity of ${\mc D}'$ along $p$, we see easily
\begin{equation}
\label{eps_ineq}
\eps\ge\mu-\nu_0.
\end{equation}
Then, since $L=k$, the strict transform $\tilde B_k$ of the line $B_k$ intersects the divisor $\tilde E$ at a point that is different from $\tilde p\cap\tilde E$. Let $h$ be a general element of the linear system $|s+2f|$ on $\tilde E$ that passes through both the points $\tilde p\cap\tilde E$ and $\tilde B_k\cap\tilde E$, where $s$ and $f$ are the minimal section and a fiber of $\tilde E$. Taking into account (\ref{DE_eq}), for a general $D\in{\mc D}$ we obtain
$$
\tilde D|_{\tilde E}\circ h=(k+2)\mu+n-\nu_1-\ldots-\nu_{k+1}-\nu_0-\eps.
$$
Denote $\theta=\frac{1}{k+1}(\nu_1+\ldots+\nu_{k+1})$. Then, if $\nu_0\ge\mu$, we get
$$
\theta\le\frac{k+2}{k+1}n.
$$
We have the same estimation even if $\nu_0<\mu$, using $\eps\ge\mu-\nu_0$. Anyway, since $k\ge 2$, we find $\theta\le\frac43n$, and from the N\"other-Fano inequality (\ref{NF_ineq}) we see that $N>3L$. Combining this estimation with (\ref{kvad1_ineq}), one gets
$$
L\left(\frac{k-1}{k}\alpha+\frac{6n^2}{k}\right)>\frac{16}3n^2L,
$$
and we have a contradiction with the estimation (\ref{alpha2_ineq}). The general infinitely near case is dealt.

\medskip

\noindent{\bf The special infinitely near case.} In this case $L=1$, and we can always assume $N\le 3$. Indeed, if $N\ge 4$, we immediately get a contradiction with the quadratic inequality (\ref{kvad1_ineq}).

The case can be dealt exactly in the same way as in \cite{IP}, so here we only give a couple of remarks. Consider a unique plane $P\subset\PA^5$ that contains the point $B_0$ and the infinitely near line $B_1$. If $\xi\not\in P$, we do not need any changes with respect to \cite{IP}. Suppose $\xi\in P$, and let $M$ be a general hyperplane in $\PA^5$ that contains $P$. Then the surface $H=M\cap V$ is a K3 surface with a unique singular (double) point at $\xi$. We follow the original explanations in \cite{IP}, but first we blow up the singular point of $H$. It is not difficult to observe that we get even more strong estimations. The unique case that can not be dealt is the case of three different lines $P\cap V$. Two of them must have their common point at $\xi$. But this situation is prohibited by the condition ${\mc G}_3$.

\medskip

\noindent{\bf The "point" case.} Let the linear system ${\mc D}$ has the multiplicity $\nu_1>2n$ at the point $B_0\ne\xi$, and there is a line $l\subset V$ that contains $B_0$ and $\xi$. Our consideration follows to the ones in \cite{IP}, with simplifications.

Let $T$ be a three-linear subspace in $\PA^5$ that is tangent to $V$ at $B_0$. From the condition ${\mc G}_4$ it follows that $Q|_T$ is a non-degenerate quadratic cone with the vertex at $B_0$. The restriction $V|_T$ consists of the lines $l=l_1,\ldots,l_k$ that are generators of the cone, and the residual curve $C$. The important observation is that the scheme-theoretic intersection $T\cap V$ is reduced along any of the lines $l_1,\ldots,l_k$ by the condition ${\mc G}_2$. Notice also that $k\le 3$ by the condition ${\mc G}_5$. The residual curve $C$ is reduced and irreducible (if we blow up the vertex of the cone $Q|_T$, it becomes a section of the corresponding ruled surface of type $\FA_2$).

Denote $\mu_i=\mult_{l_i}{\mc D}$, $i=1,\ldots,k$, and $\alpha=\mult_C{\mc D}$. Consider a general hyperplane $M\supset T$, its restriction $H=M\cap V$, and a general element $D\in{\mc D}$. Then:
\begin{itemize}
\item $H$ is a K3 surface with double (du Val) points at $B_0$ and $\xi$;
\item $D|_H=\alpha C+\mu_1 l_1+\ldots+\mu_k l_k+R$, where $\Supp R$ does not contain any of the curves $C,l_1,\ldots,l_k$.
\end{itemize}
The last fact is important: the multiplicities of these curves in ${\mc D}|_H$ coincide with the multiplicities of the linear system ${\mc D}$ along the corresponding curves. It follows from the indicated important observation.

In fact, we can even assume that $H$ has only ordinary double points at $B_0$ and $\xi$. Indeed,
since the intersection $T\cap V$ is reduced along $l_1$ and $C$, this curve is a normally 
crossing divisor at a neighborhood of $\xi$, and it remains to take into account that $V$ itself
has an ordinary double point there. Moreover, since $Q|_T$ is a non-degenerate quadratic cone 
by the condition ${\mc G}_4$, the general hyperplane section $H$ has also an ordinary double 
point at $B_0$.

Let $V'\to V$ be the blow up of the point $B_0$, and $H'$ the strict transform of $H$. 
Then $H'$ is non-singular outside of $\xi$. Mark it by $'$the strict transforms of all the 
curves $l_i$ and $C$. Denote $R'=D'|_{H'}-\alpha C'-\mu_1 l'_1-\ldots-\mu_l l'_k$. We have
$$
\begin{array}{rcl}
R'\circ C' & = & (6-k)n-(4-k)\nu_1+\frac32\alpha-\frac12\mu_1-\mu_2-\ldots-\mu_k\ge 0,\\
R'\circ l'_1 & = & n-\nu_1+\frac32\mu_1-\frac12\alpha\ge 0, \\
R'\circ l'_2 & = & n-\nu_1+2\mu_2-\alpha \ge 0,\\
& \cdots & \\
R'\circ l'_k & = & n-\nu_1+2\mu_k-\alpha \ge 0.\\
\end{array}
$$
Take the sum of all this equations,and we obtain
$$
6n+\mu_1+\ldots+\mu_k\ge (k-2)\alpha+4\nu_1>(k-2)\alpha+8n,
$$
and if $k\le 2$, we get a contradiction with the condition $\mu_i\le n$ for all $i$ (we suppose there are no maximal singularities along curves). Assume $k=3$, then we see that
$$
\frac12\mu_1+\mu_2+\mu_3\ge \frac76\alpha+\frac43\nu_1-\frac43n,
$$
and combining it with the first equation, we obtain
$$
\frac{13}3n+\frac13\alpha\ge\frac73\nu_1>\frac{14}3n,
$$
and this yields a contradiction since $\alpha\le n$.

Notice that we assume $k\le 3$ by the condition ${\mc G}_5$ (the method does not work for $k=4$, the reader can check it by himself). Lemma \ref{point_prop} is proved.

\section{Maximal singularities centered at curves.}
\label{excl2_sec}

Let $V\in{\mc F}$ be the complete intersection of a quadric $Q$ and a cubic $T$ with an ordinary double point $\xi$, as before. We suppose all the conditions ${\mc G}_1,\ldots,{\mc G}_5$ are satisfied.

\begin{remark}
\label{remark:space-curves} Let $\Lambda$ be a three-dimensional
linear subspace in $\PA^5$. Then $V\vert_{\Lambda}$
is reduced along every its non-plane component. Indeed, the
intersection $V\vert_{\Lambda}$ does not contain non-plane curves
in the case when $\Lambda\cap Q$ is reducible. On the other hand,
a cubic surface in $\PT$ does not intersect an
irreducible quadric surface by a double twisted cubic.
\end{remark}

We fix a linear system ${\mc D}\subset |-nK_V|$ that has no
fixed components.

Let $B$ an irreducible curve on $V$, and $\Lambda$ be the linear
span of a curve $B$ considered as a linear subspace~in~$\PA^5$.
Put $\nu=\mult_{B}(\mathcal{D})$.

\begin{proposition}
\label{proposition:curves} Suppose that $\nu>n$, i.e., ${\mc D}$
has a maximal singularity along $B$. Then one of the following
holds:
\begin{enumerate}
\item $B$ is a line;%
\item $B$ is a conic, and $\Lambda\subset Q$;%
\item $B$ is a conic, and $\xi\in B$.%
\end{enumerate}
\end{proposition}

\noindent{\bf Proof.} The proof of the proposition is similar to the proof of lemma 3.6 in \cite{IP} but we do our calculation on $V$ (not on its blow ups), because the claim of the proposition is simpler than the claim of lemma 3.6 in \cite{IP}, which is used not only to exclude curves on the non-singular complete intersection of a quadric and a quartic but also to find  relations between birational involutions.

We may assume that $\xi\in\Lambda$ due to \cite{IP}. Moreover, we
assume that  $\xi\in B$, because the proof is simpler in the case
when $\xi\not\in B$. The proof consists of several lemmas.

\begin{lemma}
\label{lemma:plane-cubic} The curve $B$ is not a plane cubic.
\end{lemma}

\noindent{\bf Proof.} Suppose that $B$ is a plane cubic. Then
$\Lambda\subset Q$. Let $\Lambda$ be a general a three-dimensional
subspace containing $B$. Then $\Lambda\cap V=B\cup\bar{B}$, where
$\bar{B}$ is a plane cubic.

The curves $B$ and $\bar{B}$ intersect in three distinct points,
different from $\xi$. Hence, for a general surface
$D\in\mathcal{D}$ we have
$$
3n=D\cdot\bar{B}\ge 3\nu,
$$
which is a contradiction. The lemma is proved.

We may assume that $\dim \Lambda\ge 3$. The inequality $\nu>n$
implies that $\deg B\le 5$.

\begin{lemma}
\label{lemma:smooth-big-degree} The following cases are
impossible:
\begin{itemize}
\item the curve $B$ is a rational normal curve of degree $4$ such that $\Lambda=\mathbb{P}^{4}$;%
\item the curve $B$ is a rational normal curve of degree $5$ such that $\Lambda=\mathbb{P}^{5}$;%
\item the curve $B$ is an elliptic normal curve of degree $5$ such that $\Lambda=\mathbb{P}^{4}$.%
\end{itemize}
\end{lemma}

\noindent{\bf Proof.}
Suppose $B$ is smooth. Put $d=\deg B$. Take the smallest
natural number $m\le d$ such that the following conditions
hold:
\begin{itemize}
\item the curve $B$ is cut out on the threefold $V$ in a set-theoretic sense by surfaces of the linear system $|-mK_{V}|$ that pass through the curve $B$;%
\item the scheme-theoretic intersection of two sufficiently general surfaces of the linear system $|-mK_{V}|$ passing through $B$ is reduced in the general point of $B$.%
\end{itemize}

Let $\psi\colon V'\to V$ be the extremal blow up of the curve $B$
(see \cite{Tz03}), $E$ be the exceptional divisor of $\psi$, and
$H'=\psi^{*}(-K_{V})$. Then
$$
\big(mH'-E\big)\cdot \big(\mu H'-\nu E\big)^{2}\ge 0,
$$
because the proper transform of the linear system $\mathcal{D}$ on
the threefold $V'$ does not have fixed components, but $mH'-E$ is
nef (see Lemma~5.2.5 in \cite{CPR}). Hence, the inequality
\begin{equation}
\label{equation:VA}
6mn^{2}-dm\nu^{2}-2d\nu n-n^{2}\Big(2-2g(B)-d-\frac{1}{2}\Big)\ge 0%
\end{equation}
holds (see Lemma 15 in \cite{Me03}, and the proof of Lemma~2 in \S
3 in \cite{IP}).

Putting $m=2$ in the inequality~\ref{equation:VA}, we conclude the
proof. The lemma is proved.

Thus, we see that either $\Lambda=\mathbb{P}^{3}$, or $B$ is a
curve of degree $5$ such that $\Lambda=\mathbb{P}^{4}$.

\begin{lemma}
\label{lemma:smooth-quintic} Either $\Lambda=\mathbb{P}^{3}$, or
the curve $B$ is singular.
\end{lemma}

\noindent{\bf Proof.} Suppose that the curve $B$ is a smooth
rational curve of degree $5$ such that $\Lambda=\mathbb{P}^{4}$,
which implies that $B$ is an image of a smooth rational curve of
degree $5$ in $\mathbb{P}^{3}$ via general projection. The curve
$B$ is smooth. Thus, we can use the assumptions and notations of
the proof of Lemma~\ref{lemma:smooth-big-degree}. Putting $m=3$ in
the inequality~\ref{equation:VA}, we obtain a contradiction. The
lemma is proved.

\begin{lemma}
\label{lemma:singular-quintic} Either $\Lambda=\mathbb{P}^{3}$, or
the curve $B$ is smooth at $\xi$.
\end{lemma}

\noindent{\bf Proof.} We may assume that the curve $B$ is a curve
of degree $5$ such that $\Lambda=\mathbb{P}^{4}$, and the curve
$B$ is singular at the point $\xi$. Let $W$ be a sufficiently
general surface of the linear system $|-3K_{V}|$ that contains the
curve $B$. Put $\nu_0=\mult_{\xi}\mathcal{D}$.

Let $g_0\colon \bar{V}\to V$ be a blow up of the singular point
$\xi$, $\bar{W}$ be the proper transform of the surface $W$ on
$\bar{V}$, $\bar{B}$ be a proper transform of the curve $B$ on
$\bar{V}$, and $E_0$ be the exceptional divisor of $g_0$. Then
$\deg {\mc N}_{\bar{B}|\bar{V}}=1$.

Let $g\colon\tilde{V}\to\bar{V}$ be a blow up of the curve
$\bar{B}$. Then the linear system
$$
\Big|(g_0\circ g)^*\big(-n K_{V}\big)-\nu_0 g^*(E_0)-\nu E\Big|
$$
does not have fixed components, where $E$ is the $g$-exceptional
divisor, but the complete linear system $|g^*(\bar{W})-E|$ does
not have base curves. Thus, we have
$$
\Big((g_0\circ g)\big(-n K_{V}\big)-\nu_0g^*(E_0)-\nu E\Big)^2\Big((g_0\circ g)\big(-3K_V\big)-g^*(E_0)-E\Big)\ge 0,%
$$
which implies that
$$
0\le 18n^2-10n\nu-12\nu^2+4\nu\nu_0-\nu_0^2=(18n^2-10n\nu-8\nu^2)-(2\nu-\nu_0)^2<0,%
$$
which is a contradiction. The lemma is proved.

\begin{lemma}
\label{lemma:three-dimensional-case} The curve $B$ is not a curve
of degree $5$ that is smooth at $\xi$ such that
$\Lambda=\mathbb{P}^{4}$.
\end{lemma}

\noindent{\bf Proof.} Suppose that the curve $B$ is a singular
rational curve of degree $5$ that is singular at some point
$p_1\in B$ such that $p_1\neq \xi$ and $\Lambda=\mathbb{P}^{4}$.
Let us use the arguments of the proof of
Lemma~\ref{lemma:singular-quintic} to derive a
contradiction\footnote{In fact, we can apply here use the
arguments of the proof of Lemma~\ref{lemma:smooth-big-degree}. The
inequality~\ref{proposition:curves} holds in the case when $B$ is
smooth at $\xi$, but instead of genus $g(B)$ we must plug in the
arithmetic genus of the curve $B$, which gives a contradiction,
because the curve $B$ is cut out by cubic hypersurfaces.}.

Let $g_0:\bar{V}\to V$ be a blow up of the points $\xi$ and
$p_{1}$, $\bar{W}$ be the proper transform of the surface $W$ on
$\bar{V}$, $\bar{B}$ be a proper transform of the curve $B$ on
$\bar{V}$, and $E_0$ and $E_1$ be the exceptional divisors over the points $\xi$ and $p_1$. Then $\deg {\mc N}_{\bar{B}|\bar{V}}=-2$.

Let $g:\tilde{V}\to\bar{V}$ be a blow up of $\bar{B}$. Then the
linear system
$$
\Big|(g_0\circ g)^*\big(-n K_{V}\big)-\nu_0g^*(E_0)-\nu_1g^*(E_1)-\nu E\Big|%
$$
has no fixed components, where  $E$ is the $g$-exceptional
divisor, $\nu_0=\mult_{\xi}\mathcal{D}$ and $\nu_1=\mult_{p_1}{\mc D}$.

The linear system $|g^*(\bar{W})-E|$ does not have base curves.
Thus, taking the intersection index of a general element of this linear system with two general elements of the system above, we have
$$
\begin{array}{r}
0\le 18n^2-10n\nu-14\nu^2+2\nu\nu_0-\nu_0^2+4\nu\nu_1-\nu_1^2=\\
=(18n^2-10n\nu-9\nu^2)-(\nu-\nu_0)^2-(2\nu-\nu_1)^2<0,%
\end{array}
$$
which is a contradiction. The lemma is proved.

Thus, we proved that $\Lambda=\mathbb{P}^{3}$. Put
$\Lambda=\Lambda$. Then
$$
V\vert_{\Lambda}=B\cup \bar{B},
$$
where $\bar{B}$ is a curve on the complete intersection $V$ such
that $\bar{B}\ne B$ (see Remark~\ref{remark:space-curves}). The
quadric $Q\vert_{\Lambda}$ is either a smooth quadric surface, or
an irreducible quadric cone, because $\Lambda$ is the linear span
of the irreducible curve $B$.

\begin{remark}
\label{remark:twisted-cubic} We consider only the case when $B$ is
a smooth rational curve of degree $3$. The other cases are much simpler and left to the reader.
\end{remark}

Let $\Upsilon$ be the hyperplane in $\PA^5$, that is tangent to
the quadric $Q$ at the point $\xi$, $H$ be a sufficiently general
hyperplane section of the threefold $V$ that passes through the
curve $B$, and $\pi\colon S\to H$ be the minimal resolution of
singularities of the surface $H$.

\begin{lemma}
\label{lemma:singularity-at-p0} The point $\xi$ is a singular
point of type $\mathbb{A}_{k}$ on the surface $H$. The inequality
$k\leq 2$ holds if $\Lambda\not\subset\Upsilon$.
\end{lemma}

\noindent{\bf Proof.} Let $\zeta\colon U\to V$ be a blow up of the
point $\xi$, $E$ be the exceptional divisor of $\zeta$, and
$\tilde{H}$ be the proper transform of the surface $H$ on the
threefold $U$. Then $E\cong\mathbb{P}^{1}\times\mathbb{P}^{1}$,
and the birational morphism $\pi\colon S\to H$ can be factorized
through an induced morphism
$\zeta\vert_{\tilde{H}}\colon\tilde{H}\to H$. Either the
intersection $E\cap\tilde{H}$ is irreducible, or the intersection
$E\cap\tilde{H}$ consists of two fibers of two different
projections $E\to\mathbb{P}^{1}$, respectively. In the former
case, we have $k=1$, but in the latter case we see that $\xi$ is
an isolated Du Val singular point of the surface $H$ of type
$\mathbb{A}_{k}$, and $\zeta\vert_{\tilde{H}}$ is a partial
resolution of singularities of the surface $H$.

Suppose that $\Lambda\not\subset\Upsilon$. Let us show that $k\leq
2$. We may assume that $k\ne 1$, which implies that
$E\cap\tilde{H}=E_{1}\cup E_{2}$. Put $O=E_{1}\cap E_{2}$. To
conclude the proof, we must show that the surface $\tilde{H}$ is
smooth at the point $O$.

Let $[y_0:y_1:\ldots:y_5]$ be the coordinates in $\PA^5$. Then $V$
can be defined by the equation
$$
\left\{\aligned
&y_5h(y_0:y_1:y_2:y_3:y_4)-q_1(y_0:y_1:y_2:y_3:y_4)=0,\\
&y_5q_2(y_0:y_1:y_2:y_3:y_4)-t(y_0:y_1:y_2:y_3:y_4)=0,\\
\endaligned
\right.
$$
where $h$, $q_{i}$ and $t$ are homogeneous polynomial of degree
$1$, $2$ and $3$, respectively. The point $\xi$ is given by
$y_1=y_2=y_3=y_4=0$, the quadric $Q$ is given by the equation
$y_5h-q_1=0$, the cubic $T$ is given by the equation $y_5q_2-t=0$,
and $\Upsilon$ is given by the equation $h=0$.

Let $\phi\colon V\dasharrow X$ be a projection from $\xi$, where
$X$ is a quartic threefold in $\mathbb{P}^{4}$. Then $X$ can be
naturally defined by the equation
$$
ht-q_1q_2=0,
$$
where $[y_0:y_1:y_2:y_3:y_4]$ are coordinates on $\PA^5$. The
threefold $X$ has $12$ different ordinary double points given by
the equalities $h=t=q_1=q_2=0$, $\phi(\Upsilon)$ is a hyperplane
$h=0$, and $\phi(\Lambda)$ is a plane that is not contained in
$\phi(\Upsilon)$. There is a commutative diagram
$$
\xymatrix{
&U\ar@{->}[ld]_{\zeta}\ar@{->}[rd]^{\eta}&\\%
V\ar@{-->}[rr]_{\phi}&&X,}
$$
where $\eta$ is a birational morphism that contracts $12$ smooth
rational curves to $12$ singular points of $X$, respectively. Then
$\eta(E)$ is a quadric surface on $X$ that is given by the
equations $h=q_1=0$, the curves $\eta(E_{1})$ and $\eta(E_{1})$
are lines on $\eta(E)$ such that
$\eta(E_{1})\cap\eta(E_{2})=\eta(O)$. Moreover, the plane
$\phi(\Lambda)$ intersect the quadric $\eta(E)$ either by the line
$E_{1}$, or by the line $E_{2}$.

We may assume that $\phi(\Lambda)\cap\eta(E)=E_{1}$. Put
$\breve{H}=\phi(H)$. Then $\phi(H)$ is a general hyperplane
section of the threefold $X$ that contains $X\cap\phi(\Lambda)$.
To prove that $\tilde{H}$ is smooth at the point $O$, it is enough
to prove that $\phi(H)$ is smooth at $\eta(O)$, which follows from
the fact that the one-cycle $X\cdot\phi(\Lambda)$ is reduced along
$E_{1}$. On the other hand, the cycle $X\cdot\phi(\Lambda)$ is not
reduced along $E_{1}$ if and only if the cubic $t=0$ contains the
line $E_{1}$, which implies that the line $E_{1}$ contains
two-points singular points of $X$ that are cut out on $E_{1}$ by
the equation $q_{2}=0$. However, it follows from conditions ${\mc
G}_2$ and ${\mc G}_3$ that the line $E_{1}$ contains at most one
singular point of $X$. The lemma is proved.

Note, that $\mathrm{mult}_{\xi}(B\cup\bar{B})\geq 4$ if
$\Lambda\subset\Upsilon$.

\begin{lemma}
\label{lemma:triple-line} The cycle $\bar{B}$ is not a triple
line.
\end{lemma}

\noindent{\bf Prof.}
Suppose $\bar{B}=3L$, where $L$ is a line on $H$. Put
$\bar{Q}=Q\vert_{\Lambda}$ and $\bar{T}=T\vert_{\Lambda}$. Then
$$
\bar{Q}\cdot\bar{T}=3L+B,
$$
and $\bar{Q}$ is irreducible.

Suppose that $\bar{Q}$ is smooth. Then
$\bar{Q}\cong\mathbb{P}^{1}\times\mathbb{P}^{1}$, and $B$ must be
a divisor of type $(3,0)$ on the quadric $\bar{Q}$, which is
impossible, because $B$ is irreducible and reduced.

The quadric $\bar{Q}$ is a cone, and $L$ is on of its rulings.
Then either the cubic $\bar{T}$ is singular along $L$, or the
cubic $\bar{T}$ is tangent to the quadric $\bar{Q}$ along $L$.
Hence, there is a two-di\-men\-si\-o\-nal linear subspace
$\Omega\subset\mathbb{P}^{5}$ that is tangent to both $\bar{T}$
and $\bar{Q}$ along  $L$.

The sub-scheme $V\vert_{\Omega}$ is not reduced along $L$, which
contradicts the condition ${\mc G}_2$. The lemma is proved.

It follows from Lemma~1 in \S 4 in \cite{IP} that $\Lambda$ is a
tangent linear subspace to $V$ in at most two points outside of
the singular point $\xi$. In fact, the subspace $\Lambda$ is a
tangent linear subspace to $V$ in at most one point, because
otherwise the quadric $Q\vert_{\Lambda}$ is reducible.

\begin{remark}
\label{remark:degenerate-point} Suppose that $\Lambda$ is a
tangent linear subspace to $V$ at a point $O\ne\xi$. Then $O$ is
an isolated ordinary double point of  $H$, and
$\mathrm{mult}_{O}(B\cup\bar{B})\geq 4$.
\end{remark}

\begin{corollary}
\label{corollary:triple-line} The surface $H$ has at most isolated
ordinary double points outside of the singular point $\xi$.
\end{corollary}

The birational morphism $\pi$ contracts a chain or smooth rational
curves $C_{1},\ldots, C_{k}$ to the point $\xi$ such that
$C_{i}^{2}=-2$ on the surface $S$. We may assume that
$$
C_{1}\cdot C_{2}=C_{k}\cdot C_{k-1}=1
$$
and $C_{1}\cdot C_{i}=C_{k}\cdot C_{j}=0$ for $i\ne 2$ and $j\ne
k-1$.

\begin{lemma}
\label{lemma:reducible-and-reduced} The curve $\bar{B}$ is
reducible.
\end{lemma}

\noindent{\bf Proof.}
Suppose that $\bar{B}$ is irreducible. Then $\bar{B}$ is a smooth
rational cubic curve. Let us restrict our linear system
$\mathcal{D}$ to the surface $H$. We have
$$
\mathcal{D}\vert_{H}=\nu
B+\mult_{\bar{B}}(\mathcal{D})\bar{B}+\mathcal{B}\equiv n\Big(B+\bar{B}\Big),%
$$
where $\mathcal{B}$ is a linear system on $H$ that has no fixed
components. Thus, we have
$$
0\le \Big(\big(\nu-\mu\big) B+\mathcal{B}\Big)\cdot \bar{B}=\Big(n-\mult_{\bar{B}}(\mathcal{D})\Big)\bar{B}^{2},%
$$
which implies that $\mult_{\bar{B}}(\mathcal{D})>n$ in
the case when $\bar{B}^{2}<0$. It is easy to check that the
inequality $\mult_{\bar{B}}(\mathcal{D})>n$ is
impossible, because
$$
(\nu-n) B+\mathcal{B}=\Big(n-\mult_{\bar{B}}(\mathcal{D})\Big)\bar{B}%
$$
on the surface $H$. The point $\xi$ is an intersection point of
$B$ and $\bar{B}$, and the linear subspace $\Lambda$ is not a
tangent subspace to $V$, because multiplicity of $\bar{B}\cup B$
in every point is at most two. Similarly, we see that
$\Lambda\not\subset\Upsilon$. Thus, we have $\bar{B}^{2}<0$,
because the $\bar{B}$ on $H$ can be contracted to a Du Val
singular point of type $\mathbb{A}_{k+1}$. The lemma is proved.

It follows from the proof of
Lemma~\ref{lemma:reducible-and-reduced} that to conclude the proof
it is enough to prove that the intersection form of the
irreducible components of the curve $\bar{B}$ on the surface $H$
is negatively defined. The negative definiteness of the
intersection form of the irreducible components of the curve
$\bar{B}$ implies the existence of a commutative diagram
$$
\xymatrix{
\bar{V}\ar@{-->}[rr]^{\beta}\ar@{->}[d]_{\alpha}&&\check{V}\ar@{->}[d]^{\gamma}\\%
V\ar@{-->}[rr]_{\psi}&&\mathbb{P}^{1},}
$$
where $\psi$ is a rational map induced by the projection from
$\Lambda$, $\alpha$ is an extremal blow up of the curve $B$,
$\beta$ is a composition of anti-flips in the irreducible
components of the curve $\bar{B}$, and $\gamma$ is a fibration
into $K3$ surfaces. Thus, the geometrical naturs of the negative
definiteness of the intersection form of the irreducible
components of the curve $\bar{B}$ is the existence of so-called
\emph{bad link} in the notations of \cite{CPR}.

\begin{lemma}
\label{lemma:non-reduced} The curve $\bar{B}$ is not reduced.
\end{lemma}

\noindent{\bf Proof.}
Suppose that $\bar{B}$ is reduced. Then $\bar{B}$ is either a
union of three lines, or a union of a line and a conic. We
consider only the former case. Thus, we have
$$
\bar{B}=L_{1}\cup L_{2}\cup L_{3},
$$
where $L_{1}$, $L_{2}$ and $L_{3}$ are different lines.

The lines $L_{1}$, $L_{2}$, $L_{3}$ do not pass through a smooth
point of the surface $H$, because otherwise they must lie on a
plane in $\mathbb{P}^{5}$, which contradicts to
$\Lambda=\mathbb{P}^{3}$. Thus, there are three possible subcases:
\begin{itemize}
\item the lines $L_{1}$, $L_{2}$, $L_{3}$ pass through the point $\xi$;%
\item the lines $L_{1}$, $L_{2}$, $L_{3}$ pass through a singular point of $H$ that is different from $\xi$;%
\item the lines $L_{1}$, $L_{2}$, $L_{3}$ do not intersect in one point.%
\end{itemize}

To conclude the proof it is enough to prove that the intersection
form of the irreducible components of the lines $L_{1}$, $L_{2}$
and $L_{3}$ on the surface $H$ is negatively defined.

Suppose that $\Lambda\not\subset\Upsilon$.  Then
$\Lambda\cap\Upsilon$ is a plane, which implies that at least one
line among $L_{1}$, $L_{2}$ and $L_{3}$ does not pass through the
point $\xi$ by the condition ${\mc G}_2$. Similarly, we see that
it follows form the conditions ${\mc G}_2$ and ${\mc G}_3$ that at
most one line among $L_{1}$, $L_{2}$ and $L_{3}$ passes through
the point $\xi$ (see the proof of
Lemma~\ref{lemma:singularity-at-p0}.

Let $\bar{L}_{i}$ be the proper transform of the line $L_{i}$ on
the surface $S$. Then the curves
$\bar{L}_{1},\bar{L}_{2},\bar{L}_{3},C_{1},\ldots,C_{k}$ form a Du
Val graph of type $\mathbb{D}_{k+4}$ in the case when the lines
$L_{1}$, $L_{2}$, $L_{3}$ pass through a singular point of the
surface $H$ that is different from the point $\xi$, otherwise they
form a Du Val graph of type $\mathbb{D}_{k+3}$. Hence, the lines
$L_{1}$, $L_{2}$ and $L_{3}$ can be contracted to a Du Val
singular point, which implies that they intersection form on the
surface $H$ is negatively defined (see \cite{Ar62}).

We may assume that $\Lambda\subset\Upsilon$. Then
$Q\vert_{\Lambda}$ is a cone, whose vertex is $\xi$, which implies
that $H$ is smooth outside of $\xi$, and $L_{1}$, $L_{2}$, $L_{3}$
pass through the point $\xi$.

In the case $k=1$, the lines $L_{1}$, $L_{2}$, $L_{3}$ can be
contracted to a Du Val singular point of type $\mathbb{D}_{4}$,
which implies that they intersection form on the surface $H$ is
negatively defined. Thus, we assume that $k\ne 1$. Then there is a
plane $\Omega\subset\Lambda$ such that $\Omega$ contains at least
two of the lines $L_{1}$, $L_{2}$, $L_{3}$, but $\Omega\subset Q$.
The latter is impossible by conditions conditions ${\mc G}_2$ and
${\mc G}_3$. The lemma is proved.

Thus, we may assume that $\bar{B}$ is not reduced. Then
$$
\bar{B}=2L+L^{\prime},
$$
where $L$ and $L^{\prime}$ are different lines. Moreover, it
follows from the proof of Lemma~\ref{lemma:non-reduced} that to
obtain a contradiction and conclude the proof of
Theorem~\ref{proposition:curves} it is enough to prove the
existence of a birational morphism $H\to\bar{H}$ that contracts
both lines $L$ and $L^{\prime}$.

\begin{lemma}
\label{lemma:no-degenerate-point} The linear subspace $\Lambda$ is
not a tangent linear subspace to $V$ at any smooth point of $V$.
\end{lemma}

\noindent{\bf Prof.} Suppose that $\Lambda$ is a tangent linear
subspace to $V$ at a point of $O\in V$ such that $O\ne\xi$. Then
$\mathrm{mult}_{O}(B+\bar{B})\geq 4$, which implies that $O$ is
the intersection point of the curve $L$, $L^{\prime}$ and $B$. On
the other hand, the quadric $Q\vert_{\Lambda}$ must be singular at
$O$, which implies that $Q\vert_{\Lambda}$ is an irreducible
quadric cone, whose vertex is the point $O$. Arguing as in the
proof of Lemma~\ref{lemma:triple-line}, we obtain a contradiction
with the condition ${\mc G}_2$. The lemma is proved.

Thus, the surface $H$ is smooth outside of $L\cup \xi$, and the
line $L$ passes through at most three singular points of the
surface $H$, because $H$ has isolated singularities. Moreover, the
proof of Lemma~\ref{lemma:no-degenerate-point} implies that the
quadric surface $Q\vert_{\Lambda}$ is smooth, which implies
$\Lambda\not\subset\Upsilon$.

\begin{lemma}
\label{lemma:at-most-two-points} The surface $H$ has at most two
singular points on the line $L$.
\end{lemma}

\noindent{\bf Proof.} Suppose that the surface $H$ has exactly
three singular points on $L$. Let $p_{1}$ and $p_{2}$ be singular
points of the surface $H$ contained in $L$ that are different from
the singular point $\xi$, and $\Xi$ be a hyperplane in
$\mathbb{P}^{5}$ such that $H=V\cap\Xi$. Put
$\check{Q}=Q\vert_{\Xi}$ and $\check{T}=T\vert_{\Xi}$.

Let $\Gamma$ be any hyperplane in $\Xi\cong\mathbb{P}^{4}$ that is
tangent to the quadric $\check{Q}$ at any point of the line $L$.
Then the cycle $V\vert_{\Gamma}$ is not reduced along $L$, because
$V\vert_{\Gamma}$ is a complete intersection of a quadric
$\check{Q}$ and a cubic $\check{T}$ that has at least $4$ singular
point on the line $L$.

Suppose that $\check{Q}$ is smooth. Then $\check{Q}\vert_{\Gamma}$
is a quadric cone. Hence, the arguments of the proof of
Lemma~\ref{lemma:triple-line} imply that there is a
two-dimensional linear subspace $\Omega\subset\mathbb{P}^{5}$ such
that the subschema $V\vert_{\Omega}$ is not reduced along $L$,
which is impossible by the condition ${\mc G}_2$.

Suppose that $\check{Q}$ is a cone that  is singular at some point
of $L$. Then $Q\vert_{\Lambda}$ is an irreducible cone. Arguing as
in the proof of Lemma~\ref{lemma:triple-line}, we obtain a
contradiction with the condition ${\mc G}_2$.

The quadric $\check{Q}$ must be a cone that is smooth along the
line $L$. Then $\check{Q}\vert_{\Gamma}=\Omega_{1}\cup\Omega_{2}$,
where $\Omega_{i}$ is a two-dimensional linear subspace in
$\Xi\cong\mathbb{P}^{4}$. We may assume that
$$
\Omega_{1}\cap\Omega_{2}\ne L\subset\Omega_{1},
$$
but the cycle $V\vert_{\Gamma}$ is not reduced along $L$, which
implies that the subschema $V\vert_{\Omega_{1}}$ is not reduced
along $L$, which is impossible by the condition ${\mc G}_2$. The lemma is proved.

We may assume that the surface $H$ has exactly two singular points
on the line $L$, because the other case is simpler.

\begin{lemma}
\label{lemma:three-singular-points} The line $L$ contains the
point $\xi$.
\end{lemma}

\noindent{\bf Proof.} Suppose that $\xi\not\in L$. Then $\xi\in
L^{\prime}$. Arguing as in the proof of
Lemma~\ref{lemma:non-reduced}, we see that  the lines $L$ and
$L^{\prime}$ can be contracted on the surface $H$ to a singular
point of type $\mathbb{D}_{k+4}$. The lemma is proved.

Let $p_{1}$ be the singular point of $H$ that is different from
$\xi$, and $Z$ be the curve on the surface $S$ such that
$\pi(Z)=p_{1}$. Then  $Z$ is a smooth rational curve such that
$Z^{2}=-2$, but $p_{1}\not\in L^{\prime}$, because $\Lambda$ is
not a tangent linear subspace to $V$ at $L\cap L^{\prime}$.

Let $\bar{L}$ and $\bar{L}^{\prime}$ be the proper transforms of
the curves $L$ and $L^{\prime}$ on the surface $S$, respectively.
To conclude to proof of Proposition~\ref{proposition:curves}, it
is enough to show that the curves
$$
\bar{L}, \bar{L}^{\prime}, Z, C_{1},\ldots, C_{k}%
$$
form a graph of Du Val type. The curves $\bar{L},
\bar{L}^{\prime}, Z, C_{1},\ldots, C_{k}$ form the following
graphs:
\begin{itemize}
\item a graph of type $\mathbb{A}_{4}$ in the case when $k=1$ and $\xi\in L^{\prime}$;%
\item a graph of type $\mathbb{D}_{k+3}$ in the case when $\xi\not\in L^{\prime}$.%
\end{itemize}

Thus, we may assume that $\xi\in L^{\prime}$ and $k\geq 2$. Then
the plane $\Upsilon\cap\Lambda$ must be contained in the quadric
$Q$. On the other hand, the plane $\Upsilon\cap\Lambda$ contains
both lines $L$ and $L^{\prime}$, which is impossible by conditions
conditions ${\mc G}_2$ and ${\mc G}_3$. Thus, the claim of
Proposition~\ref{proposition:curves} is proved.

\section{Remarks on points.}
\label{excl3_sec}

Let $V\in{\mc F}$ be the complete intersection of a quadric $Q$
and a cubic $T$ with an ordinary double point $\xi$, as before. We
suppose all the conditions ${\mc G}_1,\ldots,{\mc G}_5$ are
satisfied.

We fix a linear system ${\mc D}\subset |-nK_V|$ that has no fixed
components. Let $O$ be a smooth point of $V$, and $\Lambda$ be the
three-dimensional  linear subspace~in~$\PA^5$ that is tangent to
$V$ at the point $O$. Put $\nu=\mult_{O}(\mathcal{D})$.

The proof of Proposition~\ref{proposition:curves} can be adjusted
to prove the following result.

\begin{proposition}
\label{proposition:point} The inequality $\nu\leq n$ holds.
\end{proposition}

\noindent{\bf Proof.} We may assume that $\xi\in\Lambda$ due to
\cite{IP}. The quadric $Q\vert_{\Lambda}$ and the cubic
$T\vert_{\Lambda}$ are singular at $O$. Moreover, the cubic $T$ is
singular at $\xi$. Let $L_{1}$ be the line in $Q\vert_{\Lambda}$
that passes through the points $O$ and $\xi$. Then $L_{1}\subset
T$, which implies that $L_{1}\subset\Lambda\cap V$.

Let $\Upsilon$ be the hyperplane in $\PA^5$, that is tangent to
the quadric $Q$ at the point $\xi$, $H$ be a sufficiently general
hyperplane section of the threefold $V$ that contains $\Lambda\cap
V$. Then it follows from the condition ${\mc G}_5$ that
$\Lambda\not\subset\Upsilon$, which implies (see
Lemma~\ref{lemma:singularity-at-p0}) that the point $\xi$ is
either a singular point of type $\mathbb{A}_{1}$ on the surface
$H$, or a singular point of type $\mathbb{A}_{2}$ on the surface
$H$.

\begin{lemma}
\label{lemma:ODP-point-ODP} The point $\xi$ is an isolated
ordinary double point of $H$.
\end{lemma}

\noindent{\bf Proof.} Suppose that $k=2$. Then
$\mathrm{mult}_{\xi}(V\vert_{\Lambda})\geq 3$, but
$\mathrm{mult}_{O}(V\vert_{\Lambda})\geq 4$, and the cycle
$V\vert_{\Lambda}$ is reduced along $L_{1}$ by the condition ${\mc
G}_4$. Blowing up the threefold $V$ at the point $\xi$, we see
that there is a plane $\Omega\subset Q$ that contains the line
$L_{1}$ together with another line that is passed through $\xi$
and contained in $V$, which is impossible by the conditions ${\mc
G}_2$ and ${\mc G}_3$. The lemma is proved.

Let $\alpha\colon\bar{V}\to V$ be a blow up of the point $O$, and
$F$ be the exceptional divisor of $\alpha$.

\begin{lemma}
\label{lemma:ODP-point} The point $O$ is an isolated ordinary
double point of $H$.
\end{lemma}

\noindent{\bf Proof.} The claim is clear in the case when
$Q\vert_{\Lambda}$ is an irreducible quadric cone. In general, we
have the equality
$$
\mathrm{mult}_{O}(V\vert_{\Lambda})=\mathrm{mult}_{O}(H)+\sum_{Z\subset
F}\mathrm{deg}(Z)\mathrm{mult}_{Z}(\bar{H}\cdot\bar{H}^{\prime}),
$$
where $\bar{H}$ is a proper transform of $H$ on $\bar{V}$,
$\bar{H}^{\prime}$ is a proper transform on $\bar{V}$ of a general
hyperplane section of $V$ containing $V\cap\Lambda$, and $Z$ is an
irreducible curve on $F\cong\mathbb{P}^{2}$. However, it follows
from the condition ${\mc G}_2$ and ${\mc G}_3$ that the equality
$\mathrm{mult}_{O}(V\vert_{\Lambda})=4$ holds. Therefore, the
pencil spanned by the surfaces $\bar{H}$ and $\bar{H}^{\prime}$
does not has base curves contained in $F$, which implies that the
point $O$ is an isolated ordinary double point of $H$. The lemma
is proved.

It follows from the proof of Lemma~\ref{lemma:triple-line} that
$V\vert_{\Lambda}$ is reduced except possibly the case when
$Q\vert{\Lambda}$ is a union of two planes intersecting in a line
different from $L_{1}$. However, it follows from the conditions
${\mc G}_2$ and ${\mc G}_3$ that the quadric $Q\vert{\Lambda}$
does not contain a plane that contain a line different from
$L_{1}$.

\begin{corollary}
\label{corollary:reduced} The cycle $V\vert_{\Lambda}$ is reduced.
\end{corollary}

Suppose that
$$
V\vert_{\Lambda}=L_{1}+L_{1}+L_{3}+B,
$$
where $L_{i}$ is a line, and $B$ is a smooth rational curve of
degree $3$. Then $\xi\in B$, $O\in B$ and $O\in L_{i}$, because
$\mathrm{mult}_{\xi}(V\vert_{\Lambda})\geq 2$ and
$\mathrm{mult}_{O}(V\vert_{\Lambda})\geq 4$.

\begin{remark}
\label{remark:just-two-singular-points} The surface $H$ is smooth
outside of $O$ and $\xi$.
\end{remark}

Let $\pi\colon S\to H$ be the minimal resolution of singularities
of the surface $H$. Then $\pi$ is induced by a composition of
$\alpha$ and a blow of $\xi$. The birational morphism $\pi$
contracts two smooth rational curves $E$ and $Z$ such that
$\pi(E)=\xi$ and $\pi(Z)=O$. Let $\bar{L}_{i}$ be the proper
transform of the line $L_{i}$ on the surface $S$, and $\bar{B}$ be
the proper transform of the curve $B$ on the surface $S$. Then
$$
\bar{L}_{1}^{2}=\bar{L}_{2}^{2}=\bar{L}_{3}^{2}=\bar{B}^{2}=E^{2}=Z^{2}=-2,
$$
because $S$ is a smooth $K3$ surface.

It follows from the proof of Proposition~\ref{proposition:curves}
that to prove the inequality $\nu\le n$ it is enough to prove that
the intersection form of the curves $\bar{L}_{1}$, $\bar{L}_{2}$,
$\bar{L}_{3}$, $\bar{B}$ and $E$ on the surface $S$ is negatively
defined. In fact, the inequality $\nu\le n$ follows from the
semi-negative definiteness of the intersection form of the curves
$\bar{L}_{1}$, $\bar{L}_{2}$, $\bar{L}_{3}$, $\bar{B}$ and $E$ on
the surface $S$.

\begin{remark}
\label{remark:semi-negative} The semi-negative definiteness of the
intersection form of the curves $\bar{L}_{1}$, $\bar{L}_{2}$,
$\bar{L}_{3}$, $\bar{B}$ and $E$ implies the existence of a
commutative diagram
$$
\xymatrix{
\bar{V}\ar@{-->}[rr]^{\beta}\ar@{->}[d]_{\alpha}&&\check{V}\ar@{->}[d]^{\gamma}\\%
V\ar@{-->}[rr]_{\psi}&&\mathbb{P}^{1},}
$$
where $\psi$ is a rational map induced by the projection from
$\Lambda$, $\beta$ is a composition of anti-flips, and $\gamma$ is
either an elliptic fibration, or a fibration into $K3$ surfaces.
Thus, the semi-negative definiteness of the intersection form of
the curves $\bar{L}_{1}$, $\bar{L}_{2}$, $\bar{L}_{3}$, $\bar{B}$
and $E$ implies the existence of so-called \emph{bad link} in the
notations of \cite{CPR}.
\end{remark}

The curves $\bar{L}_{1}$, $\bar{L}_{2}$, $\bar{L}_{3}$, $\bar{B}$,
and $E$ form a graph $\mathbb{D}_{5}$. Therefore,  the
intersection form of the curves $\bar{L}_{1}$, $\bar{L}_{2}$,
$\bar{L}_{3}$, $\bar{B}$ and $E_{1}$ is negatively defined, which
concludes the proof of Proposition~\ref{proposition:point} in the
case when the intersection $V\cap\Lambda$ consists of three
different lines and a smooth rational cubic curve. In general, the
proof is similar and left to the reader.

\section{Finishing the proof.}
\label{fin_sec}

Let $V\in{\mc F}$ be the complete intersection of a quadric $Q$ and a cubic $T$ in $\PA^5$ with an ordinary double point $\xi\in V$. We suppose $V$ satisfies the conditions ${\mc G}_1,\ldots,{\mc G}_5$.

\begin{lemma}
\label{bir_lem}
Let $B\subset V$ be either a line, or a conic such that $\langle B\rangle\subset Q$, or a conic such that $\langle B\rangle\not\subset Q$ and $\xi\in B$. Consider a linear system ${\mc D}\subset|n(-K_X)-\nu B|$ without fixed components. Then there is a birational automorphism $\tau:V\dasharrow V$ such that:
\begin{itemize}
\item $\tau_*^{-1}{\mc D}\subset|(4n-3\nu)(-K_V)-(5n-4\nu)B|$, if $B$ is a line;
\item $\tau_*^{-1}{\mc D}\subset|(13n-12\nu)(-K_V)-(14n-13\nu)B|$, if $B$ is a conic and $\langle B\rangle\subset Q$;
\item $\tau_*^{-1}{\mc D}\subset|(15n-14\nu)(-K_V)-(16n-15\nu)B|$, if $B$ is a conic, $\langle B\rangle\not\subset Q$, and $\xi\in B$.
\end{itemize}
\end{lemma}

\noindent{\bf Proof.} The first two cases are the same as in the non-singular case, even if $B$ contains the singular point, and they are dealt completely in \cite{IP}.

Let $B$ be a conic, $\xi\in B$, and $\langle B\rangle\not\subset Q$. Then we define $\tau$ as follows. Consider a general 3-dimensional linear subspace $L\supset\langle B\rangle$. Then $V|_L=B+C$, where $C$ is an elliptic curve that intersects $B$ at the point $\xi$ and at another 3 points. Let $V'\to V$ be the blow-up of $\xi$ and then the strict transform of the curve $B$ with the exceptional divisors $E$ and $F$ respectively. We see that there is a Zariski open subset $\bar{V'}\subset V'$ that is fibred into elliptic curves over an open subset $\bar E\in E$ (really we have to throw out only a finitely many points in $E$). The there is a reflection on $\bar{V'}$ that is on each fiber $s\subset{\bar{V'}}$ a usual reflection on elliptic curve with respect to the point $s\cap\bar E$. Notice that $F$ becomes a three-section of the fibration. On $V$ this reflection gives the birational automorphism $\tau$. The reader can check that the action of $\tau$ is described as above. The
  lemma is proved.

\begin{lemma}
\label{birtransf_lem}
Let $V_1$ and $V_2$ be varieties from ${\mc F}$ that are birational to each other, and $\psi:V_1\dasharrow V_2$ is the corresponding birational map (it is described in section \ref{main_sec}). Let ${\mc D}\subset|n(-K_{V_1})-\nu\xi|$ be a linear system without fixed components. Then
$$
\psi_*{\mc D}\subset|(2n-\nu)(-K_{V_2})-(3n-2\nu)\xi|.
$$
\end{lemma}

\noindent{\bf Proof.} It can be checked immediately using the construction of $\psi$. Notice that on $V_1$ there is a unique element in $|-K_{V_1}|$ that has the multiplicity 3 at the point $\xi$. The lemma is proved.

\medskip

\noindent{\bf Proof of the part {\it (ii)} of theorem \ref{main_th}.} Let $\rho:U\to S$ be a Mori fibration, $V_1$ is a variety from ${\mc F}$ with its birational "counterpart" $V_2$, and $\chi:V_1\dasharrow U$ a birational map. Consider a very ample linear system
$$
{\mc D}_U=|n'(-K_U)+\rho^*(A)|,
$$
where $A$ is an ample divisor on $S$. Denote ${\mc D}$ the strict transform on $V_1$ of ${\mc D}_U$ by means of $\chi$.

It only remains to collect the data of propositions
\ref{point_prop} and \ref{proposition:curves} and lemmas
\ref{bir_lem} and \ref{birtransf_lem}. If ${\mc D}$ has no maximal
singularities, then $\chi$ must be an isomorphism by theorem 4.2
in \cite{Corti}. If ${\mc D}$ has maximal singularities, we apply
finitely many times the birational automorphisms from lemma
\ref{bir_lem} or the birational map from lemma \ref{birtransf_lem}
in any order. Each time we decrease the degree of the strict
transforms of ${\mc D}$, and finally we get the result by
propositions \ref{proposition:curves} and \ref{point_prop}. The
theorem is proved.

\section{Final remarks.}
\label{rem_sec}

Birational rigidity may fall under small deformations not only for Fano varieties. We can give at least two examples for del Pezzo fibrations.

The first example. Let $Q\subset\PQ$ be a non-degenerate quadratic cone, and $Q_S$ its non-singular section by a general quartic hypersurface $F$. Consider a double cover $V'\to Q$ branched along $Q_S$. The variety $V'$ is not $\QA$-factorial and has two double points over the vertex of the cone. There are two small resolutions $V_1$ and $V_2$ of $V'$ that are related by means of flop $\psi:V_1\dasharrow V_2$. It is easy to see that these varieties are actually fibrations $V_1/\POn$ and $V_2/\POn$ into del Pezzo surfaces of degree 2, and the flop $\psi$ is not square. The fibers of these fibrations arise from two pencils of planes in $Q$. It is proved (\cite{Grin1}) that, in general case, these fibrations are unique Mori fibrations in their class of birational equivalence, up to square birational maps. It is easy to see that for a general choice of $F$ the varieties $V_1$ and $V_2$ are not isomorphic to each other, thus they are non-rigid. But they become isomorphic if the se
 ction $Q_S$ is invariant with respect to the natural involution on $Q$ that exchange the pencils of planes, and in this case, $V_1/\POn$ (or, which is the same, $V_@/\POn$) is birationally rigid. Clearly, small deformations may break the symmetry of $Q_S$, and we obtain the general, i.e., non-rigid, case.

The second example. This is the case $(\eps,n_1,n_2,n_3)=(0,2,2,2)$ in \cite{Grin2}. Consider a $\PT$-fibration $X\to\POn$ that is defined by the projectivization of the bundle ${\mc E}={\mc O}\oplus{\mc O}(2)\oplus{\mc O}(2)\oplus{\mc O}(2)$. Let $Q\subset X$ be a threefold that is fibred over $\POn$ into non-degenerate quadratic cones. The vertices of the cones lie on the minimal section $t$ of $X$ (i.e., that corresponds to surjection ${\mc E}\to{\mc O}$), $Q\sim 2M-4L$, where $M$ is the tautological divisor, $L$ is a fiber of $X$. Consider a double cover $V\to X$ that is branched over a non-singular section $Q_R=Q\cap R$, where $R\sim 3M$. Notice that $t\circ R=0$. We see that $V$ is a fibration $V/\POn$ into del Pezzo surfaces of degree 1. It has a section $s$ that lies over $t$, and there is a flop $\psi:V\dasharrow U$ centered at $s$, onto a fibration $U/\POn$ into del Pezzo surfaces of degree 1 with the same construction. It is proved that, up to square birational map
 s, these fibrations are unique Mori fibrations in their class of birational equivalence. As before, they are not isomorphic and non-rigid in general case, but for specially chosen divisors $R$, they becomes isomorphic, and thus rigid by definition \ref{rig_def}.

\end{document}